\documentclass[runningheads,a4paper]{elsarticle}
\usepackage{amsmath,amsfonts,amsthm}
\usepackage{amssymb,amscd}
\usepackage{graphicx}
\usepackage{textcomp} 

\def\muntz{M\"untz }
\def\muntzsecond{M\"untz}
\def\divdif{\mathord\kern.43em\vrule width.6pt height5.6pt depth-.28pt \kern-.43em\Delta}

\newtheorem{theorem}{Theorem}

\newtheorem{proposition}{Proposition}
\newtheorem{definition}{Definition}
\newtheorem{corollary}{Corollary}

\newdefinition{remark}{Remark}

\begin{document}

\title{Dimension Elevation in \muntz Spaces:
A New Emergence of the \muntz Condition}

\author[]{Rachid Ait-Haddou \corref{cor}}
\ead{Rachid.AitHaddou@kaust.edu.sa}

\cortext[]{Corresponding author}

\address{King Abdullah University of Science and Technology, Thuwal, Saudi Arabia}

\begin{abstract}
We show that the limiting polygon generated by the dimension elevation
algorithm with respect to the \muntz space $span(1,t^{r_1},t^{r_2},...,t^{r_m},...)$,
with $0 < r_1 < r_2 < ... < r_m < ...$ and $\lim_{n\to\infty}r_n = \infty$, over an interval
$[a,b]\subset]0,\infty[$ converges to the underlying Chebyshev-B\'ezier curve if and only if
the \muntz condition $\sum_{i=1}^{\infty} \frac{1}{r_i} = \infty$ is satisfied.
The surprising emergence of the \muntz condition in the problem raises
the question of a possible connection between the density questions
of nested Chebyshev spaces  and the convergence of the corresponding
dimension elevation algorithms. The question of convergence with no condition of
monotonicity or positivity on the pairwise distinct real numbers $r_i$ remains an open problem.
\end{abstract}

\maketitle
\section{Introduction}
Degree elevation of B\'ezier curves is a standard
technique in computer aided curve design. It consists of
iteratively expressing a B\'ezier curve of a fixed degree in the
Bernstein bases of the linear spaces of polynomials of
higher degrees. The process generates a sequence of control
polygons which converges uniformly to the underlying B\'ezier curve
\cite{Prau}. Degree elevation, or more appropriately,
dimension elevation, can be generalized to any infinite nested sequence
of linear spaces in which an analogue notion of Bernstein basis
can be defined. More precisely, let $E_{\infty} =
(u_1,u_2,...,u_m,...)$
be an infinite sequence of sufficiently differentiable functions
$u_i$ over an interval $[a,b]$ and such that for every
$n \geq 1$, the linear space
$DE_{n} = span(u'_1,u'_2,...,u'_n)$ is an extended
Chebyshev space of dimension $n$ over $[a,b]$ \cite{mazure1}. Then,
for any $n \geq 1$, the linear space $E_{n} = span(1,u_1,...,u_n)$
possesses a so-called Chebyshev-Bernstein basis $B^n_k,
k=0,...,n$ and characterized as the unique normalized basis
of $E_{n}$ such that for every $k \in \{1,..,n\}$, the
function $B^n_{k}$ has $k$ zeros at $a$ and $(n-k)$ zeros
at $b$ \cite{mazure2}. A Chebyshev-B\'ezier curve, $\Gamma$,
in the space $E_{n}$ can be expressed in the Chebyshev-Bernstein
bases of the spaces $E_{n}$ and $E_{n+1}$ over an interval $[a,b]$ as
\small
\begin{equation*}
P(t) = \sum_{i=0}^{n} B^n_i(t) P_i =
\sum_{i=0}^{n+1} B^{n+1}_i(t) P^{1}_i, \quad P_i,  P^{1}_i \in \mathbb{R}^s; \quad s\geq 1.
\end{equation*}
\normalsize
The defining endpoint conditions of Chebyshev-Bernstein
bases show that the points $P^{1}_i$ are related to the points $P_i$ as
follows: $P^{1}_0 = P_{0}$, $P^{1}_{n+1} = P_{n}$ and for
$i=1,...,n$, there exist real numbers $\xi_i\in ]0,1[$ such
that
\begin{equation}\label{affinity}
P^{1}_i = (1 -\xi_i) P_{i-1} + \xi_i P_{i}.
\end{equation}
Iterating the process of expressing the curve $\Gamma$  in the
Chebyshev-Bernstein bases of the nested sequence of spaces
$E_{n+1} \subset  E_{n+2} \subset  ... \subset E_{m} \subset ...$
generates a sequence of control polygons. Since dimension elevation is a corner cutting scheme,
the generated sequence converges to a Lipschitz-continuous curve \cite{deboor}. However,
a difficult question is to characterize the sequences $E_{\infty}$ in which the limiting polygon
converges uniformly to the underlying Chebyshev-B\'ezier curve.
Here, we report our complete solution to the problem when $E_{\infty}$ is the \muntz sequence
$E_{\infty} = (t^{r_1},t^{r_2},...,t^{r_n},...)$ with $0 < r_1 < r_2 < ... < r_m < ...$ 
and $\lim_{n\to\infty}r_n = \infty$. More precisely, we have 
\begin{theorem}\label{maintheorem}
The limiting polygon generated by the dimension elevation
algorithm with respect to the \muntz space $span(1,t^{r_1},t^{r_2},...,
t^{r_m},...)$ with $0 < r_1 < r_2 < ... < r_m < ...$ and $\lim_{n\to\infty}r_n = \infty$
over an interval $[a,b]\subset]0,\infty[$ converges uniformly to the
underlying non-constant Chebyshev-B\'ezier curve if and only if
\begin{equation}\label{muntzcondition}
\sum_{i=1}^{\infty} \frac{1}{r_i} = \infty.
\end{equation}
\end{theorem}

Let us compare our theorem with the celebrated original \muntz Theorem
on the density of \muntz spaces \cite{almira,gurariy,muntz}
\begin{theorem} {\bf{(\muntz Theorem)}}\label{muntztheorem}
Let $( r_1, r_2,..., r_m, ...)$ be an infinite strictly increasing sequence
of positive real numbers such that $\lim_{n\to\infty} r_{n} = \infty$.
The \muntz space $span(1,t^{r_1},...,t^{r_m},...)$ is a dense
subset of $C([0,1])$ (the linear space of continuous functions on $[0,1]$
endowed with the uniform norm) if and only if
\begin{equation*}
\sum_{i=1}^{\infty} \frac{1}{r_i} = \infty.
\end{equation*}
\end{theorem}

The emergence of the \muntz condition (\ref{muntzcondition}) in both Theorem
\ref{maintheorem} and Theorem \ref{muntztheorem} is rather
surprising and may suggest a deep connection between the problem of
density of nested Chebyshev spaces and the convergence of the associated
dimension elevation algorithms. For nested \muntz spaces over the interval $[0,1]$, a hypothesis of equivalence 
is ruled out by the fact that the condition $\lim_{n\to\infty} r_n = \infty$ can be dropped
in \muntz Theorem \ref{muntztheorem} \cite{borwein}, while such condition is necessary
for the convergence of the dimension elevation over $[0,1]$ to the underlying curve 
(See Theorem \ref{muntzelevation}). However, the necessity of the condition $\lim_{n\to\infty} r_n = \infty$  
in Theorem \ref{maintheorem} remains an open problem.

It is interesting to note that as much as the classical Weierstrass approximation Theorem is a special
case of \muntz Theorem for the exponents $r_n = n$, the classical convergence theorem of the degree
elevation of B\'ezier curves is a consequence of Theorem \ref{maintheorem} for
the same set of exponents.

The \muntz condition (\ref{muntzcondition}) also appears in different contexts other than the density questions
of \muntz spaces, such as  in Biernack Theorem on entire functions \cite{marden}, in Ramm-Horv\'ath Theorem
on the inverse scattering problem \cite{horvath} or in the famous Erd\H{o}s conjecture on arithmetic
progressions \cite{erdos}.

\bigskip

\noindent {\bf{The strategy for the proof of Theorem \ref{maintheorem}: }}For nested \muntz spaces, 
the analytical form of $\xi_i$ in (\ref{affinity})
can be expressed in terms of a quotient of generalized Schur functions that depend on the interval
parameters $a$ and $b$ \cite{aithaddou1}. Iterating the dimension elevation process
using these quotients leads to complicated expressions that hinder a direct proof of
Theorem \ref{maintheorem}. The limiting curve generated by the dimension elevation in
\muntz spaces over an interval $[a,b]$ depends only on the shape parameter $b/a$
\cite{mazure1}. Therefore, we only need to prove Theorem \ref{maintheorem} over an
interval $[a,1]$ with $0 < a < 1$. Taking into account the special role of the origin in the
density questions of \muntz spaces \cite{borwein}, we could first look at the dimension
elevation algorithm in \muntz spaces over the interval $[0,1]$. An apparent obstruction to
such strategy is the fact that \muntz spaces does not possess Chebyshev-Bernstein bases over the
interval $[0,1]$. Fortunately, as we will show in this work, the pointwise limits of the
Chebyshev-Bernstein bases over an interval $[a,1]$ as $a$ goes to zero do exist and in fact
coincide with the Gelfond-Bernstein bases of \muntz spaces \cite{aithaddou2}. The dimension elevation algorithm
for the Gelfond-Bernstein bases over the interval $[0,1]$ can be easily expressed.
This allows us to prove Theorem \ref{maintheorem} over the interval $[0,1]$.
To prove Theorem \ref{maintheorem} over an interval $[a,1]$ such that $0 < a < 1$, we use the notion
of Chebyshev blossoming in \muntz spaces in order to show that the control polygon
of the dimension elevation over $[a,1]$ can be obtained by a generalized
de Casteljau algorithm from the control polygon of the dimension elevation over the interval $[0,1]$.
The necessity of the \muntz condition is proved with the aid of the bounded Chebyshev inequality
in \muntz spaces. Although the initial steps of the proof of Theorem 1, namely the proof of the theorem over the
interval $[0,1]$, appeared in our work in \cite{aithaddou2,aithaddou3}, we review here of all the steps of the proof
in order to give a consistent presentation and show how we extended our methods to deal with the away from the origin case.       

\section{Chebyshev-Bernstein Bases in \muntz Spaces}
Throughout this work, we will denote by $\Lambda_{\infty} = (r_0=0,r_1,...,r_m,...)$
an infinite sequence of strictly increasing real numbers.
For any integer $n$, we denote by $\Lambda_{n}$ the
finite subsequence $\Lambda_{n} = (r_0=0,r_1,...,r_n)$ and
denote by $E(\Lambda_n) = span(t^{r_0}=1,t^{r_1},...,
t^{r_n})$ the associated \muntz space. To give an explicit expression
for the Chebyshev-Bernstein basis of the linear space $E(\Lambda_n)$
over an interval $[a,b]$, we introduce the following definitions and terminology
\begin{definition}
A finite sequence of real numbers $\lambda = (\lambda_1,\lambda_2,...,\lambda_n)$  is termed 
a {\it{real partition}} if it satisfies
\begin{equation*}
\lambda_1 > \lambda_2 - 1 > \lambda_3 - 2 > ...> \lambda_n - (n-1) > -n.
\end{equation*}
The generalized Schur function indexed by a real partition $\lambda$ is defined as
the continuous extension of the function defined for pairwise distinct real values
$u_1,u_2,...,u_n$ by \cite{mazurejat}
\begin{equation*}\label{schurdeterminant}
S_{\lambda}(u_1,...,u_n) = \frac{\det ( u_i^{\lambda_j + n - j} )_{
{1 \leq i, j \leq n}}}{ \prod_{1\leq i < j \leq n}(u_i - u_j)}.
\end{equation*}
\end{definition}
We use the notation $S_{\lambda}(u_1^{m_1},u_2^{m_2},...,u_k^{m_k})$
to mean the evaluation of the generalized Schur function in which the argument $u_1$
is repeated $m_1$ times, the argument $u_2$ is repeated $m_2$ times and so on.
When the elements of the finite sequence $\lambda$ are positive integers, we recover the
classical notion of integer partitions whose associated Schur functions
$S_{\lambda}(u_1,...,u_n)$ are elements of the ring $\mathbb{Z}[u_1,...,u_n]$.
The value $S_{\lambda}(1^n)$ can be computed for integer partitions
using the hook-length formula \cite{sagan} or in general by the formula
\begin{equation*}\label{generalhook}
S_{\lambda}(1^{n}) = \frac{\prod_{1\leq j < k \leq n}
(\lambda_{j} - \lambda_{k} -j + k)}{\prod_{j=1}^{n}(j-1)!}.
\end{equation*}
\begin{definition}\label{realpartition}
For a finite sequence $\Lambda_n = (r_0,r_1,...,r_{n})$ of $(n+1)$ real
numbers such that $0 =r_0 < r_1 < ... < r_n$, we define the real
partition $\lambda = (\lambda_1,...,\lambda_{n},\lambda_{n+1})$ associated with
the finite sequence $\Lambda_n$ by
\begin{equation*}
\lambda_k = r_n - r_{k-1} - (n-k+1)
\quad \textnormal{for} \quad k=1,...,n+1.
\end{equation*}
We also denote by $\lambda^{(0)}$ the real partition
$\lambda^{(0)} = (\lambda_2,\lambda_3,...,\lambda_{n+1})$ termed
the bottom partition of $\lambda$.
\end{definition}
With the above definitions, we can prove the following theorem \cite{aithaddou1}

\begin{theorem}\label{bernsteintheorem}
The Chebyshev-Bernstein basis
$(B^{n}_{0,\Lambda_n},...,B^{n}_{n,\Lambda_n})$
over an interval $[a,b]$ of the \muntz space  $E(\Lambda_n)$  is given by
\begin{equation*}\label{bernstein}
B^{n}_{k,\Lambda_n}(t) = \frac{S_{\lambda}(1^{n+1})}{S_{\lambda^{(0)}}(1^{n})}
B^{n}_{k}(t)
\frac{S_{\lambda^{(0)}}(a^{n-k},b^{k}) t^{\lambda_1}
S_{\lambda}(a^{n-k},b^{k},\frac{ab}{t})}
{S_{\lambda}(a^{n+1-k},b^{k}) S_{\lambda}(a^{n-k}, b^{k+1})},
\end{equation*}
where $B^{n}_{k}$ is the classical Bernstein basis of the polynomial
space over the interval $[a,b]$, $\lambda$ the real partition associated with $\Lambda_n$
and $\lambda^{(0)}$ the bottom partition of $\lambda$.
\end{theorem}
\subsection{Gelfond-Bernstein Bases as Pointwise Limits of Chebyshev-Bernstein Bases}
Let $f$ be a smooth real function defined on an interval $I$.
For any real numbers $x_0 < x_1 < ...< x_n$ in the interval $I$,
the divided difference $[x_0,...,x_n]f$ of the function $f$ supported
at the point $x_i, i=0,...,n$ is recursively defined by $[x_0]f = f(x_0)$
and
\begin{equation*}\label{dvdefinition}
[x_0,x_1,...,x_n]f = \frac{[x_1,...,x_n]f - [x_0,x_1,...,x_{n-1}]f}{x_n - x_0}
\quad \textnormal{if} \quad n>0.
\end{equation*}
Consider, now, the function $f_t$ defined on $[0,+\infty[$ as
\begin{equation}\label{functionft}
\begin{cases}
 f_t(x) = t^x \quad \textnormal{for} \quad t>0,\; x\geq 0 \\
f_0(0) :=1, f_0(x)= 0 \quad \textnormal{for} \quad x >0.
\end{cases}
\end{equation}
\begin{definition}
For a finite sequence $\Lambda_n=(0=r_0,r_1,...,r_n)$ of strictly
increasing positive real numbers, the Gelfond-Bernstein basis of the
\muntz space $E(\Lambda_n)$ with respect to the interval $[0,1]$
is defined by
$$
H^{n}_{k,\Lambda_n}(t) = (-1)^{n-k} r_{k+1}...r_{n} [r_k,...,r_n]f_{t}
\quad \textnormal{for} \quad k=0,...,n-1
$$
and
$$
H^{n}_{n,\Lambda_n}(t) = t^{r_n},
$$
where $f_t$ is the function defined in (\ref{functionft}).
\end{definition}
The relation between the Chebyshev-Bernstein basis and the Gelfond-Bernstein
basis of a given \muntz space is given by the following theorem \cite{aithaddou2}
\begin{theorem}\label{maintheoremgelfond}
Let $\Lambda_n=(0=r_0,r_1,...,r_n)$ be a finite sequence of strictly
increasing real numbers. We denote by $B^{n}_{k,\Lambda_n}, k=0,...,n$ the Chebyshev-Bernstein
basis of the \muntz space $E(\Lambda_n)$ over an interval $[a,1]$. Then, for $k=0,...,n$ and $t\in [0,1]$ we have
$$
\lim_{a\to 0} B^{n}_{k,\Lambda_n}(t) = H^{n}_{k,\Lambda_n}(t).
$$
\end{theorem}
\noindent
The proof of Theorem \ref{maintheoremgelfond} is based on applying to the explicit
expression of the Chebyshev-Bernstein bases in Theorem \ref{bernsteintheorem},
the following splitting formula for generalized Schur functions:
If $\eta = (\lambda_1, . . . , \lambda_k, \mu_1, . . . ,
\mu_h)$ is a real partition then we have
\begin{equation*}\label{split}
\lim_{\epsilon\to 0} \frac{S_{\eta}(z_1,...,z_k,\epsilon y_1,
...,\epsilon y_h)} {\epsilon^{|\mu|}} =
S_{\lambda}(z_1,...,z_k) S_{\mu}(y_1,...,y_h),
\end{equation*}
where $\lambda$ and $\mu$ are the real partitions
$\lambda = (\lambda_1,...,\lambda_k)$ and $\mu = (\mu_1,...,\mu_h)$ and where
$|\mu|$ denotes $\mu_1 + \mu_2 +....+\mu_h$. The proof also relies
on the following interesting connection between generalized Schur functions and divided
differences of the function $f_t$ in (\ref{functionft}). Namely, for any finite sequence
$\Lambda_n = (0=r_{0},r_{1},...,r_{n})$ of strictly
increasing real numbers we have
$$
[r_0,r_1,...,r_n]f_{t} = \frac{(-1)^{n}}{r_1 r_2 ... r_n}
(1-t)^n \frac{S_{\lambda}(1,t^{n})}{S_{\lambda^{(0)}}(t^{n})},
$$
where $\lambda$ is the real partition associated with the finite
sequence $\Lambda_n$ and $\lambda^{(0)}$ the bottom
partition of $\lambda$.
\subsection{Dimension Elevation in \muntz Spaces over the Interval $[0,1]$}
To express the corner cutting scheme associated with the dimension
elevation of Gelfond-B\'ezier curves, we should express the Gelfond-Bernstein
basis of $E(\Lambda_n)$ in terms of the Gelfond-Bernstein basis of $E(\Lambda_{n+1})$.
Such expression is given by the following proposition \cite{aithaddou2}

\begin{proposition}\label{dimensionproposition}
For $k=0,...,n$, and for any $t \in [0,1]$ we have
\begin{equation*}
H^{n}_{k,\Lambda_n}(t) = \frac{r_{n+1} - r_{k}}{r_{n+1}}
H^{n+1}_{k,\Lambda_{n+1}}(t) + \frac{r_{k+1}}{r_{n+1}}
H^{n+1}_{k+1,\Lambda_{n+1}}(t).
\end{equation*}
\end{proposition}
Let $P$ be an element of the \muntz space $E(\Lambda_n)$ written as
\begin{equation}\label{expansion2}
P(t) = \sum_{k=0}^{n} H_{k,\Lambda_n}^{n}(t) P_{k} =
\sum_{k=0}^{m} H_{k,\Lambda_m}^{m}(t) b^{m}_{k}; \quad m > n.
\end{equation}
By Proposition \ref{dimensionproposition}, the control points $b^{m}_{k} = P^{m-n}_{k}$ in (\ref{expansion2})
can be computed using the following corner cutting scheme:
For $i=0,1,...,n$, we set $P_{i}^{0} =P_{i}$ and for
$j=1,2,...m-n$, we construct iteratively new polygons
$(P_0^{j},P_1^{j},....,P_{n+j}^j)$ using the inductive
rule
\begin{equation}\label{initial}
P_{0}^{j} =P_{0}^{j-1} \quad  P_{n+j}^{j} = P_{n+j-1}^{j-1}
\end{equation}
and for $i=1,...,n+j-1$
\begin{equation}\label{cornercutting}
P_{i}^{j} = \frac{r_{i}}{r_{n+j}} P_{i-1}^{j-1} +
\left( 1-\frac{r_{i}}{r_{n+j}} \right) P_{i}^{j-1}.
\end{equation}
\begin{figure}
\centerline{\includegraphics[width=0.55\textwidth]{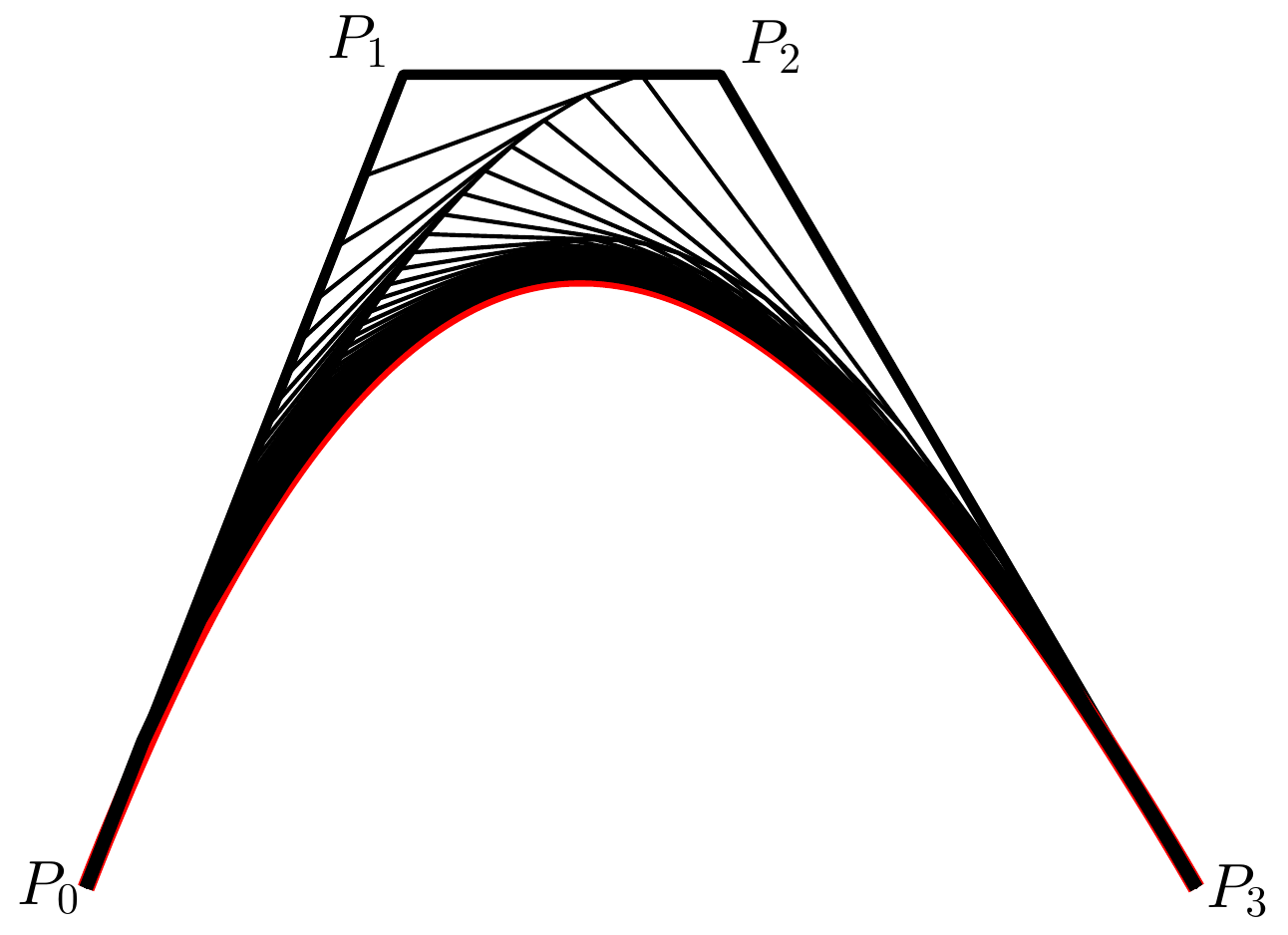}}
\caption{The sequence of control polygons generated from 100 iterations of
the corner cutting scheme (\ref{initial}) and (\ref{cornercutting})
and parameters $n=3, r_1=1, r_2=2, r_3 =3$ and $r_i = 2i$ for $i \geq 4$.
The red curve is the B\'ezier curve associated with the control polygon
$(P_0,P_1,P_2,P_3)$.}
\label{fig:figure1}
\end{figure}
When the real numbers $r_i$ are given by $r_i = i$
for every index $i$, then the corner cutting scheme
(\ref{initial}) and (\ref{cornercutting}) leads to the classical degree
elevation algorithm, in which it is well known that the limiting control
polygon converges to the underlying B\'ezier curve as $m$ goes to infinity \cite{Prau}.
Consider the case where $r_i = i$ for $i=1,...,n$ and
$r_i = 2i$ for $i > n$. Figure \ref{fig:figure1} shows the generated
polygons from 100 iterations with $n=3$. The figure clearly suggests the convergence
of the generated polygons to the B\'ezier curve with control
points $(P_0,P_1,P_2,P_3)$. Consider, now, the case where
$r_i = i$ for $i=1,...,n$, while $r_i = i^2$ for $i >n$.
Figure \ref{fig:figure2} shows the obtained polygons after 100 iterations with $n=3$.
It is clear from the figure that the limiting polygon does
not converge to the B\'ezier curve with control points $(P_0,P_1,P_2,P_3)$.
Now, consider, for example, the limiting polygon of the corner
cutting scheme (\ref{initial}) and (\ref{cornercutting})
for the case $n=3$ and in which $r_1=2, r_2 = 4, r_3 =10$ and
$r_i = 2i + 5$ for $i >3$. Figure \ref{fig:figure3} shows the generated
polygons from 100 iterations and the Gelfond-B\'ezier curve
associated with the \muntz space $F = span(1,t^{r_1},t^{r_2},t^{r_3})
= span(1,t^2,t^4,t^{10})$ and control polygon $(P_0,P_1,P_2,P_3)$.
The figure suggests that the limiting polygon converges to the Gelfond-B\'ezier curve.
To exhibit the importance of the condition $\lim_{s\to\infty} r_s = \infty$ in the dimension
elevation of Gelfond-B\'ezier curves, Figure \ref{fig:figure4} shows the limiting polygon for the case
$n = 3$, $r_1 = 1, r_2 = 2, r_3 = 3$ and
$r_i = 4 - \frac{1}{i}$ for $i>3$. The limiting polygon
does not converge to the B\'ezier curve with control polygon
$(P_0,P_1,P_2,P_3)$. In fact we can prove the following theorem \cite{aithaddou3}
\begin{figure}
\centerline{\includegraphics[width=.55\textwidth]{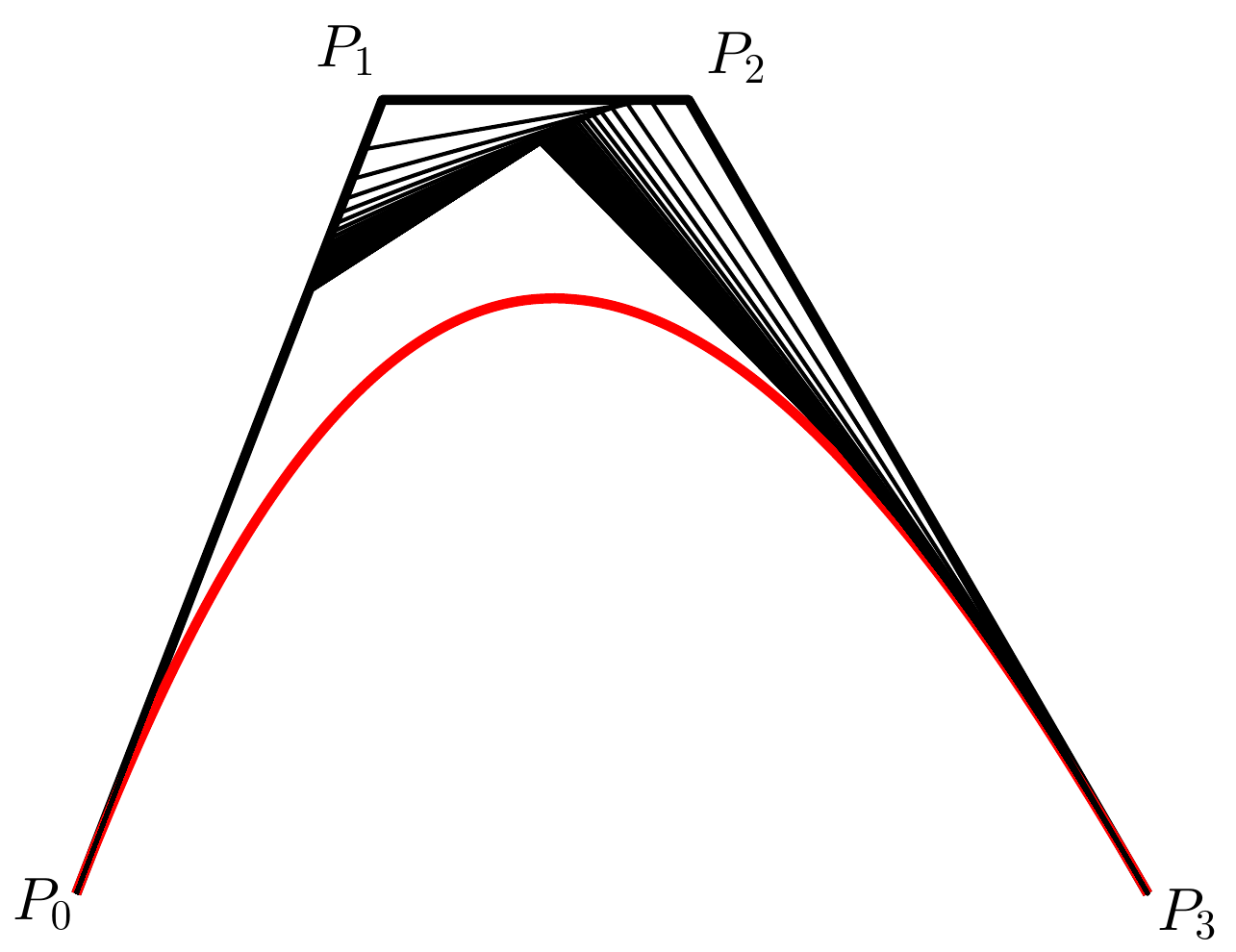}}
\caption{The sequence of polygons generated from 100 iterations of
the corner cutting scheme (\ref{initial}) and (\ref{cornercutting})
and parameters $n=3, r_1=1, r_2=2, r_3 =3$ and $r_i = i^2$
for $i \geq 4$. The red curve is the B\'ezier curve associated
with the control polygon $(P_0,P_1,P_2,P_3)$.}
\label{fig:figure2}
\end{figure}
\begin{figure}
\centerline{\includegraphics[width=.55\textwidth]{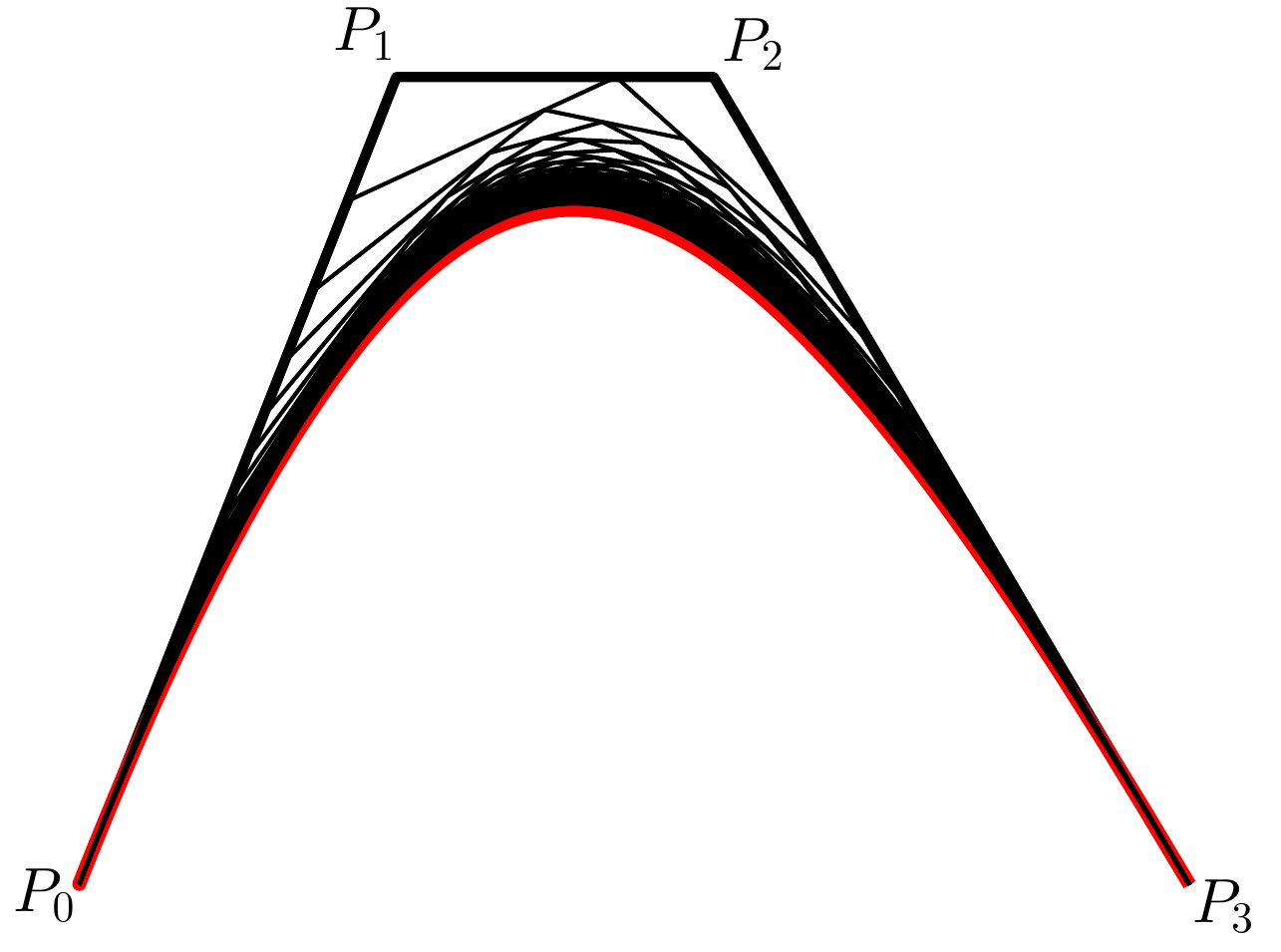}}
\caption{The sequence of polygons generated from 100 iterations of
the corner cutting scheme (\ref{initial}) and (\ref{cornercutting})
and parameters $n=3, r_1=2, r_2=4, r_3 =10$ and $r_i = 2i+5$ for $i \geq 4$.
The red curve is the Gelfond-B\'ezier curve associated with the \muntz space
$span(1,t^2,t^4,t^{10})$ and control polygon $(P_0,P_1,P_2,P_3).$}
\label{fig:figure3}
\end{figure}
\begin{figure}
\centerline{\includegraphics[width=.55\textwidth]{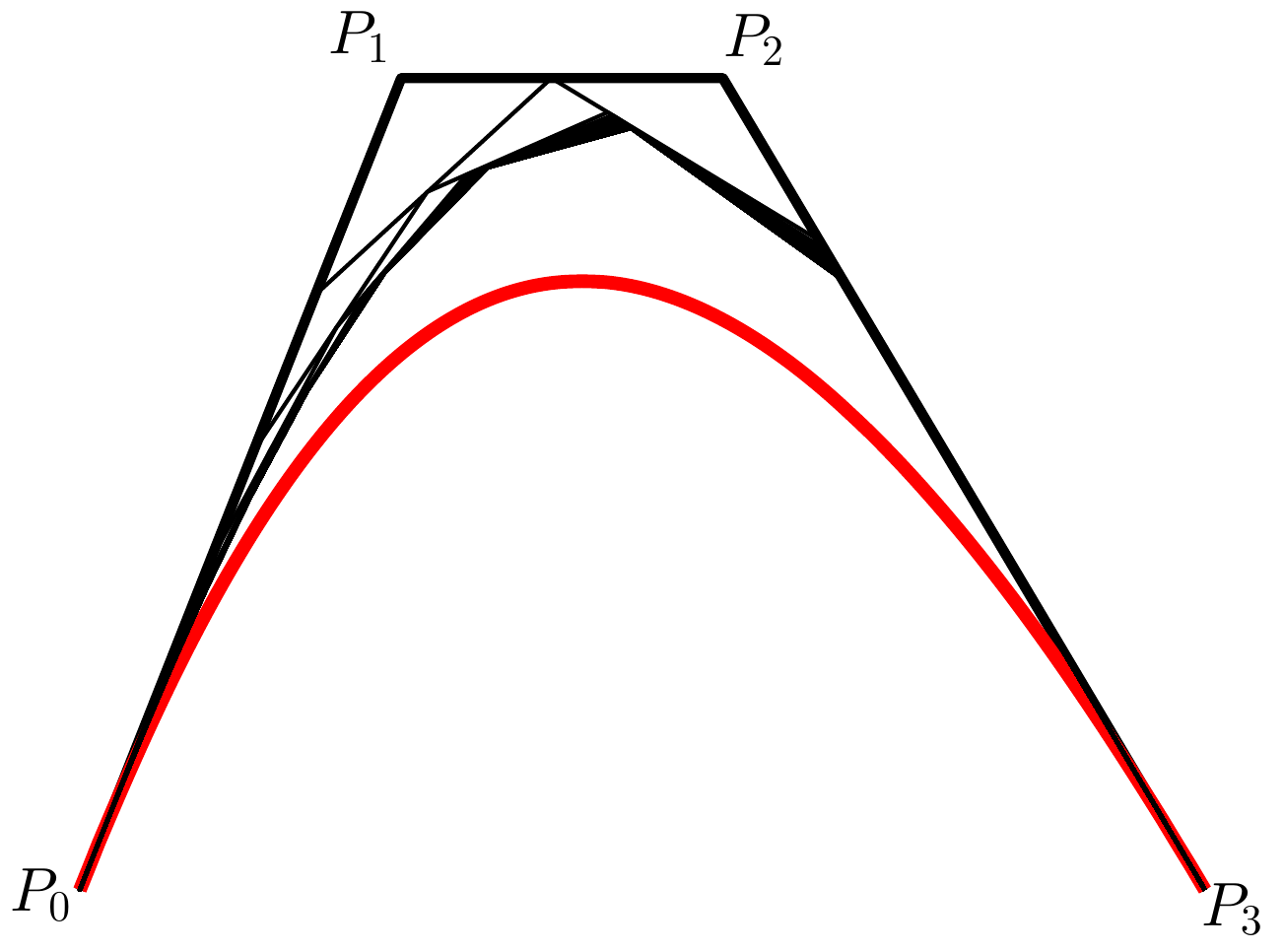}}
\caption{The sequence of control polygons generated from 100 iterations of
the corner cutting scheme (\ref{initial}) and (\ref{cornercutting})
and parameters $n=3, r_1=1, r_2=2, r_3 =3$ and $r_i = 4-\frac{1}{i}$ for $i \geq 4$.
The red curve is the B\'ezier curve associated with
the control polygon $(P_0,P_1,P_2,P_3)$.}
\label{fig:figure4}
\end{figure}
\begin{theorem}\label{muntzelevation}
The limiting polygon generated from a non-constant polygon
$(P_0,P_1,$ $...,P_n)$ in $\mathbb{R}^s, s \geq 1$ using the corner cutting scheme
(\ref{initial}) and (\ref{cornercutting}) with respect to the sequence $\Lambda_{\infty}$
converges uniformly to the Gelfond-B\'ezier curve associated with the \muntz space
$E(\Lambda_n)$ and control polygon $(P_0,P_1,...,P_{n})$ if and only if the real
numbers $r_i$ satisfy the conditions
\begin{equation*}
\lim_{s\to\infty} r_{s} = \infty \quad \textnormal{and} \quad
\sum_{i=1}^{\infty} \frac{1}{r_i} = \infty.
\end{equation*}
\end{theorem}
For the rest of the paper, we adopt the following notation: for a given element $P$ of the
\muntz space $E(\Lambda_n)$, we denote by $\eta_i(P,\Lambda_n,[0,1])$ the Gelfond-B\'ezier
control points of $P$ with respect to the interval $[0,1]$ and we denote by
$\eta_i(P,\Lambda_n,[a,1])$ the Chebyshev-B\'ezier control points of $P$ with respect
to the interval $[a,1]$. One of the main tool in proving Theorem \ref{muntzelevation},
and which will be needed later, is the fact that for any element $P$ of $E(\Lambda_n)$ we have
\begin{equation}\label{convergenceterm}
\lim_{m\to\infty} || P(\eta_{i}(t^{r_1},\Lambda_m,[0,1])^{1/r_1}) - \eta_i(P,\Lambda_m,[0,1]) ||_{\infty}  = 0
\end{equation}
uniformly in $i$.
\section{Dimension elevation over an interval $[a,1]$}
Given two intervals $[a,1] \subset [b,1], (b > 0)$, we can infer from the de Casteljau algorithm
of Chebyshev-B\'ezier curves that for any element $P$ of $E(\Lambda_n)$ there exist
real numbers $s^{(i)}_{j} \in [0,1]; j=0,...,n$ (independent of $P$) such that
$\sum_{j=0}^{n} s_{j}^{(i)} = 1$ and
\begin{equation}\label{decasteljaupositive}
\eta_{i}(P,\Lambda_n,[a,1]) = \sum_{j=0}^{n} s^{(i)}_{j}   \eta_{j}(P,\Lambda_n,[b,1]).
\end{equation}
Although the latter statement remains true for $b = 0$, its proof is not obvious and
requires the generalization of the notion of blossoming in \muntz spaces over the interval $[0,1]$
and the use of the relation between the Gelfond-Bernstein bases and the Chebyshev-Bernstein bases as stated in
Theorem \ref{maintheoremgelfond}. We have  \cite{aithaddou2}
\begin{theorem}\label{decasteljau}
Let $P$ be an element of the \muntz space $E(\Lambda_n)$. Then for any $i=0,1,...,n$
there exist real numbers $s^{(i)}_{j} \in [0,1]; j=0,...,n$ (independent of $P$) such that
\begin{equation*}
\eta_{i}(P,\Lambda_n,[a,1]) = \sum_{j=0}^{n} s^{(i)}_{j}   \eta_{j}(P,\Lambda_n,[0,1])
\; \textnormal{with} \;
\sum_{j=0}^{n} s_{j}^{(i)} = 1.
\end{equation*}
\end{theorem}
Using Theorem \ref{decasteljau}, we can prove the following theorem.
\begin{theorem}\label{uniform1}
Under the \muntz condition (\ref{muntzcondition}) on the sequence $\Lambda_{\infty}$, for any positive integer $k \leq m$ 
\begin{equation*}
\lim_{m\to\infty} | \eta_i(t^{r_1},\Lambda_m,[a,1])^{\frac{r_k}{r_1}} - \eta_i(t^{r_k},\Lambda_m,[a,1])| = 0
\end{equation*}
uniformly in $i$.
\end{theorem}
\begin{proof}
Consider the parametric curve $\Gamma: (t^{r_1},t^{r_k})$ over
the interval $[a,1]$. Since $r_k > r_1$, $\Gamma$ is the graph of a convex function. As a parametric curve in the \muntz
space $E(\Omega)$ with $\Omega = (0,r_1,r_k)$, it is thus located above its control polygon. This remain true when considering
$\Gamma$ as a parametric curve in $E(\Lambda_m)$ (as dimension elevation is a corner cutting scheme). Therefore, we have
for $i=0,...,m$
\begin{equation}\label{localinequality}
\eta_{i}(t^{r_1},\Lambda_m,[a,1])^{\frac{r_k}{r_1}} \geq \eta_{i}(t^{r_k},\Lambda_m,[a,1]).
\end{equation}
From Theorem \ref{decasteljau}, we have
\begin{equation*}
\begin{split}
& \eta_{i}(t^{r_1},\Lambda_m,[a,1])^{\frac{r_k}{r_1}} -  \eta_{i}(t^{r_k},\Lambda_m,[a,1]) = \\
& \left(\sum_{j=0}^m s_j^{(i)}  \eta_{j}(t^{r_1},\Lambda_m,[0,1])\right)^{\frac{r_k}{r_1}} -
\sum_{j=0}^m s_j^{(i)}  \eta_{j}(t^{r_k},\Lambda_m,[0,1]).
\end{split}
\end{equation*}
Therefore, by Jensen Inequality and Inequality (\ref{localinequality}), we have
\begin{equation*}
\begin{split}
& 0 \leq \eta_{i}(t^{r_1},\Lambda_m,[a,1])^{\frac{r_k}{r_1}} - \eta_{i}(t^{r_k},\Lambda_m,[a,1]) \leq  \\
& \sum_{j=0}^m s_j^{(i)}  ( \eta_{j}(t^{r_1},\Lambda_m,[0,1])^{\frac{r_k}{r_1}} - \eta_{j}(t^{r_k},\Lambda_m,[0,1])).
\end{split}
\end{equation*}
The proof is concluded upon using Equation (\ref{convergenceterm}).
\end{proof}
\begin{corollary}\label{uniform2}
Under the \muntz condition (\ref{muntzcondition}) on the sequence $\Lambda_{\infty}$, for any $P \in E(\Lambda_n)$ and
for any positive integer $k \leq m$, we have
\begin{equation*}
\lim_{m\to\infty} || P( \eta_i(t^{r_1},\Lambda_m,[a,1])^{\frac{1}{r_1}}) - \eta_i(P,\Lambda_m,[a,1])||_{\infty} = 0
\end{equation*}
uniformly in $i$.
\end{corollary}
\begin{proof}
Denote by $P$ an element of $E(\Lambda_n)$ given by
$P(t) = \sum_{k=0}^{n}  t^{r_k} A_k$. We have
\begin{equation*}
\begin{split}
& ||P(\eta_{i}(t^{r_1},\Lambda_m,[a,1])^{1/r_1}) - \eta_i(P,\Lambda_m,[a,1]||_{\infty}  \leq  \\
& \sum_{k=0}^n ||A_k||_{\infty}  | ( \eta_{i}(t^{r_1},\Lambda_m,[a,1])^{r_k/r_1} -  (\eta_{i}(t^{r_k},\Lambda_m,[a,1]) |.
\end{split}
\end{equation*}
We conclude the proof using Theorem \ref{uniform1}.
\end{proof}
Using the Widder-Hirschman-Gelfond Theorem \cite{almira,lorentz}, we proved in \cite{aithaddou3} that the point set
$D_{m} = \{\eta_i(t^{r_1},\Lambda_m,[0,1])^{1/r_1}, i=0,...,m \}$ form a dense subset of the interval $[0,1]$
as $m$ goes to infinity. Applying the generalized de Casteljau algorithm to the function $t^{r_1}$ to compute its control points over
the interval $[a,1]$ for its control points over the interval $[0,1]$ and using the density property of the point set $D_{m}$, we also have the 
following
\begin{theorem}\label{densitytheorem}
Under the \muntz condition (\ref{muntzcondition}), the point set  $D_{m} = \{\eta_i(t^{r_1},\Lambda_m,$ $[a,1])^{1/r_1},$ $ i=0,...,m \}$ form a
dense subset of the interval $[a,1]$ as $m$ goes to infinity.
\end{theorem}
\smallskip
We are now in a position to prove the main Theorem \ref{maintheorem} when the sequence
$\Lambda_{\infty}$ satisfies the \muntz condition (\ref{muntzcondition}).
In this case, if we denote by $P$ an element of $E(\Lambda_n)$,
we have to show that given a point $t \in [a,1]$ and a sequence of real numbers
$\eta_{i_m(t)}(t^{r_1},\Lambda_m,[a,1])^{1/r_1}$ that converges to $t$ as $m$ goes
to infinity (this is possible thanks to the density result in Theorem
\ref{densitytheorem}), the point $b^{m}_{i_m(t)} = \eta_{i_m(t)}(P,\Lambda_m,[a,1])$
converges to $P(t)$ as $m$ goes to infinity uniformly on $t$. We have  
\begin{equation*}
\begin{split}
& \max_{t} ||P(t) - b^{m}_{i_m(t)}||_{\infty} \leq
\max_t ||P(t) - P(\eta_{i_m(t)}(t^{r_1},\Lambda_m,
 [a,1])^{1/r_1})||_{\infty} + \\
& \max_{i_m(t)}  ||P(\eta_{i_m(t)}(t^{r_1},\Lambda_m,[a,1])^{1/r_1})-b^{m}_{i_m(t)}||_{\infty}.
\end{split}
\end{equation*}
The function $P$ is continuous in the compact interval $[a,1]$, thus
$$
\max_t ||P(t) - P(\eta_{i_m(t)}(t^{r_1},\Lambda_m,[a,1])^{1/r_1})||_{\infty} \rightarrow 0
\;  \textnormal{as} \;  m \rightarrow \infty,
$$
and Corollary \ref{uniform2} shows that
$$
\max_{i_m(t)} ||P(\eta_{i_m(t)}(t^{r_1},\Lambda_m,[a,1])^{1/r_1})-b^{m}_{i_m(t)}||_{\infty}
\rightarrow 0
\; \textnormal{as} \;  m \rightarrow \infty.
$$
This concludes the proof of the \textgravedbl if\textacutedbl  part of Theorem \ref{maintheorem}.

\begin{remark} Using Equation \ref{decasteljaupositive}, it can be proven that if the limiting polygon generated
by the dimension elevation with respect to a sequence $\Lambda_{\infty}$ over an interval $[a,b]$ converges
to the underlying Chebyshev-B\'ezier curve, then the limiting polygon over any interval $[c,d] \subset [a,b]$
also converges to the underlying curve.
\end{remark}
To prove that the \muntz condition (\ref{muntzcondition}) is necessary in Theorem \ref{maintheorem}, we need the following 
bounded Chebyshev inequality \cite{almira,gurariy,borwein}
\begin{theorem} {\bf{(Bounded Chebyshev's Inequality)}}
Let us assume that the sequence $\Lambda_\infty = (r_0 =0, r_1,r_2 ,...)$ satisfies $\sum_{i=1}^{\infty} 1/r_i < \infty$.
Then for any real-valued element $P \in E(\Lambda_{\infty})$ and for each $\epsilon >0$ there is a constant
$c(\Lambda_{\infty},\epsilon) >0$ depending only on $\Lambda_\infty$ and $\epsilon$
 (and not on the number of terms in $P$) and such that
\begin{equation}\label{bounded}
||P'||_{[0,1-\epsilon]} \leq c(\Lambda_{\infty},\epsilon)  ||P||_{[1-\epsilon,1]}.
\end{equation}
\end{theorem}
To prove that the \muntz condition (\ref{muntzcondition}) is necessary in Theorem \ref{maintheorem},
we proceed as follows: Let  $P$ be a non-constant element of $E(\Lambda_n)$ expressed in the Chebyshev-Bernstein
basis over an interval $[a,1]$ as
$$
P(t) = \sum_{i=0}^{m} B_{i,\Lambda_m}^{m}(t) b^{m}_{i} \quad m \geq n.
$$
Without loss of generality, we can assume that $P'(a) \not= 0$ (otherwise, we work on an interval
$[b,1]$ such that $[a,1] \subset [b,1]$ and $P'(b) \not=0$ and invoke Remark1).
We have \cite{aithaddou1}
$$
P'(a) = (B_{0,\Lambda_m}^{m})'(a) (b^{m}_{0} - b^{m}_{1}).
$$
Using the bounded Chebyshev inequality (\ref{bounded}) with $\epsilon = 1-a$, we have
$$
|(B_{0,\Lambda_m}^{m})'(a)| \leq c(\Lambda_{\infty},\epsilon)  ||B_{0,\Lambda_m}^{m}||_{[a,1]} =
 c(\Lambda_{\infty},\epsilon)
$$
for any $m \geq n$. Therefore,
$$
\lim_{m \to \infty} ||b^{m}_{0} - b^{m}_{1}||_{\infty} \geq \frac{||P'(a)||_{\infty}}{c(\Lambda_{\infty},\epsilon) } > 0.
$$
This shows that the dimension elevation leads to a limiting curve with a segment of non-zero length as part of the curve.
Thereby, the limiting polygon cannot converge to the underlying Chebyshev-B\'ezier curve.

\section{Concluding Remarks}
For a sequence of distinct real positive numbers $\Lambda_{\infty} = (r_0=0,r_1,...,r_n,...)$
(with no monotonicity condition on the $r_i$) the space $E(\Lambda_{\infty})$ is dense in $C([0,1])$
if and only if \cite{borwein}
\begin{equation}\label{degenerate}
 \sum_{k=1}^{\infty} \frac{r_k}{{r_k}^2 + 1} = \infty.
\end{equation}
Gelfond-B\'ezier curves are  too "degenerate'' at
the origin to study the dimension elevation algorithm
in case we have no condition of monotonicity on the
real numbers $r_i$. For instance, if we consider the case $n=3,$
$r_1=1, r_2 = 2, r_3 = 3$ and $r_j = 1/j$, for $j > 3$ and we start
with a control polygon $(P_0,P_1,P_2,P_3)$ then the control polygon
obtained by a dimension elevation to the order $m$ is not obtained
by a corner cutting scheme similar to (\ref{initial}) and (\ref{cornercutting}), but
instead the algorithm collapses the first $m-3$ control points
to $P_0$ while the remaining control points are given by $(P_1,P_2,P_3)$ \cite{aithaddou2}.
Therefore, an analogue formulation as in (\ref{degenerate}) for the dimension elevation in \muntz spaces
over $[0,1]$ is unlikely. However, if we consider the dimension elevation algorithm of Gelfond-B\'ezier curves
away from the origin, i.e., over an interval $[a,1]$ with $a >0$, then the
Gelfond-Bernstein basis coincides with the Chebyshev-Bernstein basis,
the degeneracy at the origin disappears and the algorithm leads
to a family of corner cutting schemes without imposing any condition
of monotonicity on the real numbers $r_i$. Unfortunately, such family
of corner cutting schemes involves rather complicated
coefficients expressed in term of generalized Schur functions \cite{aithaddou1}.
It will be interesting to find, for the away from the origin case,
conditions on the real number $r_i$ for the convergence of the
dimension elevation algorithm to the underlying curve.
In the theory of the density of \muntz spaces over an
interval $[a,b]$ with $a>0$, with no condition of positivity nor
monotonicity on the pairwise distinct real numbers $r_i$,
the corresponding \muntz space is a dense subset of $C([a,b])$
if and only if the real numbers $r_i$ satisfy the so-called
full \muntz condition \cite{almira,borwein}
\begin{equation}\label{fullmuntz}
\sum_{r_k \neq 0} \frac{1}{|r_k|} = \infty.
\end{equation}
An interesting question is then: Let $[a,b]$ be an interval with $a >0$ and let $r_i$
be a sequence of pairwise distinct real numbers without any condition of positivity or monotonicity.
Can we claim that the corresponding dimension elevation algorithm over
$[a,b]$ converges to the underlying curve if and only if the sequence $r_i$ satisfies the full \muntz condition (\ref{fullmuntz})?

\section*{Acknowledgments}
The author wishes to express his gratitude to Professor Marie-Laurence Mazure for pointing out an error in the original submission of this work. The author is also grateful to the referees for the many helpful comments that improved the presentation of this work.



\begin{thebibliography}{1}

\bibitem{almira}
J. M. Almira, \muntz type theorem I, Surveys in Approx. Theory 3 (2007) 152--194.

\bibitem{aithaddou1}
R. Ait-Haddou, Y. Sakane, T. Nomura, Chebyshev blossoming in \muntz
spaces: toward shaping with Young diagrams, J. Comp. Applied. Math. 247 (2013) 172--208.

\bibitem{aithaddou2}
R. Ait-Haddou, Y. Sakane, T. Nomura, Gelfond-B\'ezier curves,
Comp. Aided Geom. Design 30 (2013) 199--225.

\bibitem{aithaddou3}
R. Ait-Haddou, Y. Sakane, T. Nomura,
A \muntz type theorem for a family of corner cutting schemes,
Comput. Aided Geome. Design 30 (2013) 240--253.

\bibitem{borwein}
P. Borwein, T. Erd\'elyi, Polynomials and Polynomial Inequalities,
Graduate Texts in Mathematics, Springer, 1996.

\bibitem{deboor}
C. de Boor, Cutting corners always works,
Comput. Aided Geom. Design 4 (1987) 125--131.

\bibitem{erdos}
P. Erd\H{o}s, Problems in number theory and combinatorics, Proc. Sixth Manitoba Conf. on Num. Math.,
Congress Numer. XVIII, (1977) 35--58.

\bibitem{gurariy}
V. I. Gurariy, W. Lusky, Geometry of \muntz spaces, Lecture Notes in Maths, vol 1870, 2005.

\bibitem{horvath}
M. Horv\'ath, Inverse scattering with fixed energy and an inverse eigenvalue problem on the half-line,
Trans. Amer. Math. Soc. 358  (2006) 5161–-5177.

\bibitem{lorentz}
G-G. Lorentz, Bernstein Polynomials. University
of Toronto Press, Toronto, 1953.

\bibitem{marden}
M. Marden, Geometry of Polynomials, Math. Surveys no. 3, American Mathematical Society, Providence, RI, 1966.

\bibitem{mazure1}
M-L. Mazure, Blossoming: A geometrical approach, Constr. Approx. 15 (1999) 33–-68.

\bibitem{mazure2}
M-L. Mazure, Chebyshev-Bernstein bases, Comput. Aided Geom. Design 16 (1999) 649---669.

\bibitem{mazurejat}
M-L. Mazure, On Chebyshevian spline subdivision, J. Approx. Theory 143 (2006) 74–-110.

\bibitem{muntz}
Ch-H. \muntzsecond, \"Uber den Approximationsatz von
Weierstrass, Mathematische Abhandlungen in H. A. Schwarz's Festschrift,
Berlin, Springer (1914) 303--312.

\bibitem{Prau}
H. Prautzsch, K. Kobbelt, Convergence of subdivision and degree elevation,
Adv. comput. Math. 2 (1994) 143--154.

\bibitem{sagan}
B. Sagan, The symmetric Group: Representations, Combinatorial Algorithms and Symmetric Functions
(Second Edition) Springer, New York, 2001.

\end{thebibliography}
\end{document}